\documentclass[11pt]{article}
\usepackage{amssymb}
\usepackage{amsmath}
\usepackage{amsthm}
\usepackage{graphicx}
\usepackage{enumerate}
\usepackage[all,2cell]{xy} \UseAllTwocells \SilentMatrices
\reversemarginpar
\renewcommand\footnotemark{}
\numberwithin{equation}{section}
\theoremstyle{definition}

\DeclareMathOperator{\Sp}{Sp}
\DeclareMathOperator{\Mp}{Mp}

\begin{document}
\newcommand{\cb}[1]{\left\{#1\right\}}
\newcommand{\lb}[1]{\left(#1\right)}
\newcommand{\ls}[1]{\left[#1\right]}
\newcommand{\la}[1]{\left|#1\right|}
\newcommand{\lr}[1]{\left\langle#1\right\rangle}
\newcommand{\lat}[2]{\left.#1\right|_{#2}}
\newcommand{\dd}[2]{\frac{d{#1}}{d{#2}}}
\newcommand{\pp}[2]{\frac{\partial{#1}}{\partial{#2}}}
\newcommand{\dt}{\ensuremath{\left.\dd{}{t}\right|_{t=0}}}
\newcommand{\hb}{\ensuremath{\widehat{b}}}
\newcommand{\hc}{\ensuremath{\widehat{c}}}
\newcommand{\mg}{\ensuremath{\mathfrak{g}}}
\newcommand{\teta}{\ensuremath{\widetilde{\eta}}}
\newcommand{\hj}{\ensuremath{\widehat{j}}}
\newcommand{\mk}{\ensuremath{\mathfrak{k}}}
\newcommand{\mmp}{\ensuremath{\mathfrak{mp}}}
\newcommand{\hnu}{\ensuremath{\widehat{\nu}}}
\newcommand{\tp}{\ensuremath{\widetilde{\phi}}}
\newcommand{\hp}{\ensuremath{\widehat{\phi}}}
\newcommand{\hpsi}{\ensuremath{\widehat{\psi}}}
\newcommand{\tr}{\ensuremath{\widetilde{\rho}}}
\newcommand{\hr}{\ensuremath{\widehat{\rho}}}
\newcommand{\msp}{\ensuremath{\mathfrak{sp}}}
\newcommand{\mt}{\ensuremath{\mathfrak{t}}}
\newcommand{\bv}{\ensuremath{\mathbf{v}}}
\newcommand{\hxi}{\ensuremath{\widehat{\xi}}}
\newcommand{\bx}{\ensuremath{\mathbf{x}}}
\newcommand{\hx}{\ensuremath{\widehat{x}}}
\newcommand{\hxs}{\ensuremath{\widehat{x}^*}}
\newcommand{\txi}{\ensuremath{\widetilde{\xi}}}
\newcommand{\tchi}{\ensuremath{\widetilde{\chi}}}
\newcommand{\hchi}{\ensuremath{\widehat{\chi}}}
\newcommand{\by}{\ensuremath{\mathbf{y}}}
\newcommand{\hy}{\ensuremath{\widehat{y}}}
\newcommand{\hys}{\ensuremath{\widehat{y}^*}}
\newcommand{\bz}{\ensuremath{\mathbf{z}}}
\newcommand{\hz}{\ensuremath{\widehat{z}}}
\newcommand{\hzs}{\ensuremath{\widehat{z}^*}}
\newcommand{\oz}{\ensuremath{\overline{z}}}
\newcommand{\one}{\ensuremath{\mathbf{1}}}
\newcommand{\sC}{\ensuremath{\mathcal{C}}}
\newcommand{\C}{\ensuremath{\mathbb{C}}}
\newcommand{\CP}{\ensuremath{\mathbb{CP}}}
\newcommand{\hF}{\ensuremath{\widehat{F}}}
\newcommand{\hG}{\ensuremath{\widehat{G}}}
\newcommand{\sJ}{\ensuremath{\mathcal{J}}}
\newcommand{\N}{\ensuremath{\mathbb{N}}}
\newcommand{\tPsi}{\ensuremath{\widetilde{\Psi}}}
\newcommand{\hPsi}{\ensuremath{\widehat{\Psi}}}
\newcommand{\Q}{\ensuremath{\mathbb{Q}}}
\newcommand{\R}{\ensuremath{\mathbb{R}}}
\newcommand{\Z}{\ensuremath{\mathbb{Z}}}
\newcommand{\contr}{\ensuremath{\lrcorner\,}}
\newcommand{\Det}{\ensuremath{\mbox{Det}}}
\renewcommand{\Re}{\ensuremath{\mbox{Re}\,}}
\renewcommand{\Im}{\ensuremath{\mbox{Im}\,}}
\newcommand{\maps}[1]{\ensuremath{\stackrel{#1}{\longrightarrow}}}
\newcommand{\two}[2]{\ensuremath{\lb{\begin{array}{c}{#1}\\{#2}\end{array}}}}
\newcommand{\four}[4]{\ensuremath{\lb{\begin{array}{cc}{#1}&{#2}\\{#3}&{#4}\end{array}}}}

\newtheorem{theorem}{Theorem}[section]
\newtheorem{example}[theorem]{Example}
\newtheorem{definition}[theorem]{Definition}
\newtheorem{lemma}[theorem]{Lemma}
\newtheorem{proposition}[theorem]{Proposition}
\newtheorem{corollary}[theorem]{Corollary}

\title{Equivariant Metaplectic-c Prequantization of Symplectic Manifolds with Hamiltonian Torus Actions}
\author{Jennifer Vaughan \\ Department of Mathematics \\ University of Manitoba \\ Winnipeg, MB R3T 2N2, Canada \\ \texttt{jennifer.vaughan@umanitoba.ca}}

\maketitle

\begin{abstract}

This paper determines a condition that is necessary and sufficient for a metaplectic-c prequantizable symplectic manifold with an effective Hamiltonian torus action to admit an equivariant metaplectic-c prequantization.  The condition is evaluated at a fixed point of the momentum map, and is shifted from the one that is known for equivariant prequantization line bundles.

Given a metaplectic-c prequantized symplectic manifold with a Hamiltonian energy function, the author previously proposed a condition under which a regular value of the function should be considered a quantized energy level of the system.  This definition naturally generalizes to regular values of the momentum map for a Hamiltonian torus action.  We state the generalized definition for such a system, and use an equivariant metaplectic-c prequantization to determine its quantized energy levels.

\end{abstract}


\section{Introduction}

Metaplectic-c quantization was introduced by Hess \cite{h1} and further developed by Robinson and Rawnsley \cite{rr1}.  It is a generalization of the Kostant-Souriau quantization procedure with half-form correction that applies to a strictly broader class of symplectic manifolds.  The starting point for this paper is a symplectic manifold $(M,\omega)$ that admits a metaplectic-c prequantization.  

Suppose there is a Hamiltonian $G$ action on $(M,\omega)$ for some Lie group $G$.  Broadly speaking, a prequantization bundle for $(M,\omega)$ is called an equivariant prequantization if the $G$ action lifts to the prequantization bundle in a manner that preserves all of its structures.  This concept has been applied to prequantization line bundles in the context of a torus action and an arbitrary compact Lie group action \cite{ggk}.  It has also been applied to spin-c structures in the context of a circle action \cite{ckt} and a torus action \cite{gk}.  

In this paper, we assume that $(M,\omega)$ is metaplectic-c prequantizable and has an effective Hamiltonian torus action with momentum map $\Phi:M\rightarrow\mt^*$, where $\mt^*$ is the dual of the Lie algebra $\mt$ for the torus.  Section \ref{sec:HamMpc} contains our conventions for Hamiltonian torus actions, and reviews the definitions of the metaplectic-c group and a metaplectic-c prequantization.  

In Section \ref{sec:EquivarMpc}, we further assume that the torus action has at least one fixed point $z$.  We give the definition of an equivariant metaplectic-c prequantization, and we determine a condition on the value of the momentum map at $z$ that is necessary and sufficient for $(M,\omega)$ to admit an equivariant metaplectic-c prequantization.  For an equivariant prequantization line bundle, a comparable result is known \cite{ggk}:  the value $\frac{1}{h}\Phi(z)$ must be in the integer lattice of $\mt^*$, where $h=2\pi\hbar$ is Planck's constant.  The condition that we obtain is shifted from this due to the lift of the torus action to the symplectic frame bundle for $(M,\omega)$.  The statement of our equivariance condition is in Theorem \ref{thm:Equivar}.

In an earlier paper \cite{v2}, we defined a quantized energy level for the metaplectic-c prequantized system $(M,\omega,H)$, where $H\in C^{\infty}(M)$ is viewed as a Hamiltonian energy function on $(M,\omega)$.  This definition has a natural generalization to families of Poisson-commuting functions.  In particular, in Section \ref{sec:QuantE}, we apply it to the components of the momentum map $\Phi$ for the torus action.  Given an equivariant metaplectic-c prequantization for the system $(M,\omega,\Phi)$, we show that the regular values $x$ of $\Phi$ that are quantized are exactly those such that $\frac{1}{h}x$ lies in the integer lattice.  This is Theorem \ref{thm:EquivarQuantE}.  The section concludes by demonstrating that if the torus acts freely on the level set corresponding to a quantized energy level, then the symplectic reduction is metaplectic-c prequantizable.

Lastly, in Section \ref{sec:Ex}, we consider two examples.  First, we obtain the quantized energy levels for a harmonic oscillator of arbitrary dimension.  The result includes the half-shift predicted by the standard quantum mechanical calculation.  Then we consider the complex projective space $\CP^2$ with an action of the two-dimensional torus $T^2$ that is induced from a linear $T^2$ action on $\C^3$.  We determine the shift in the momentum map required to satisfy the equivariance condition, and find the quantized energy levels.  This example is notable because $\CP^2$ admits a metaplectic-c prequantization but not a metaplectic structure, meaning that quantization results for this system cannot be duplicated using Kostant-Souriau quantization with the half-form correction.


\section{Hamiltonian Torus Actions and Metaplectic-c Prequantization}\label{sec:HamMpc}

Section \ref{subsec:HamTs} sets up our notation and conventions for the torus and the momentum map.  In Section \ref{subsec:Mpc}, we summarize the basic elements of metaplectic-c prequantization.  Considerably more detail, including proofs, were given by Robinson and Rawnsley \cite{rr1}.  A similar review also appears in \cite{v2}.


\subsection{Hamiltonian torus actions}\label{subsec:HamTs}

Let $T^k$ be a $k$-dimensional torus with Lie algebra $\mt$.  Write $\tau\in T^k$ as $(\tau_1,\ldots,\tau_k)$ where each $\tau_j\in U(1)$.  Let $\cb{\xi_1,\ldots,\xi_k}$ be the standard basis for $\R^k$, and identify $\mt$ with $\R^k$ such that for any $\xi=\sum_{j=1}^{k}a_j\xi_j\in\R^k$, $$\exp(\xi)=(e^{2\pi ia_1},\ldots,e^{2\pi ia_k}).$$

Let $(M^{2n},\omega)$ be a connected symplectic manifold, where $n\geq k$.  For the remainder of the paper, we assume that $T^k$ has an effective Hamiltonian action on $M$  with momentum map $\Phi:M\rightarrow\mt^*$.  For all $\xi\in\mt$, we define the vector field $\xi_M$ on $M$ by $$\xi_M(m)=\lat{\dd{}{t}}{t=0}\exp(t\xi)\cdot m,\ \ \forall m\in M.$$  Our convention for the momentum map is $$d\Phi^\xi=\xi_M\contr\omega,\ \ \forall\xi\in\mt.$$  

For each $\xi\in\mt$, denote the flow of $\xi_M$ on $M$ by $\phi_\xi^t$.  Explicitly,  $$\phi^t_\xi(m)=\exp(t\xi)\cdot m,\ \ \forall m\in M,$$ and $\phi^t_\xi$ is a symplectomorphism for all $t$.  For any $\xi\in\mt$, the action of the element $\tau=\exp(\xi)\in T^k$ on $M$ is given by the map $\phi^1_\xi$.  In particular, if $\xi\in\Z^k\subset\mt$, then $\phi^1_\xi$ is the identity map on $M$.


\subsection{Metaplectic-c prequantization}\label{subsec:Mpc}

Fix a model $2n$-dimensional symplectic vector space $(V,\Omega)$, together with a compatible complex structure $J$ on $V$.  The symplectic group for $(V,\Omega)$ is denoted by $\Sp(V)$.  The metaplectic group is the connected double cover $\Mp(V)\maps{\sigma}\Sp(V)$.  The metaplectic-c group is defined to be $$\Mp^c(V)=\Mp(V)\times_{\Z_2}U(1).$$  

We will make use of the following two group homomorphisms on $\Mp^c(V)$.  The \textbf{projection map} $\sigma$ appears in the short exact sequence $$1\rightarrow U(1)\rightarrow\Mp^c(V)\maps{\sigma}\Sp(V)\rightarrow 1$$ and restricts to the double covering map on $\Mp(V)$.  The \textbf{determinant map} $\eta$ appears in the short exact sequence $$1\rightarrow\Mp(V)\rightarrow\Mp^c(V)\maps{\eta}U(1)\rightarrow 1$$ and acts on $\lambda\in U(1)\subset\Mp^c(V)$ by $\eta(\lambda)=\lambda^2$.  The Lie algebra  $\mmp^c(V)$ can be identified with $\msp(V)\oplus\mathfrak{u}(1)$ under $\sigma_*\oplus\frac{1}{2}\eta_*$.

For any $g\in\Sp(V)$, let $$C_g=\frac{1}{2}(g-JgJ).$$  Then $C_g$ commutes with $J$, so it is a complex linear map on $V$.  It can be shown \cite{rr1} that $\Det_\C C_g\neq 0$ for all $g\in\Sp(V)$.

We define an embedding of $\Mp^c(V)$ into $\Sp(V)\times\C$ such that each $a\in\Mp^c(V)$ is mapped to the pair $(g,\mu)\in\Sp(V)\times\C$, where $\sigma(a)=g$ and $\eta(a)=\mu^2\Det_{\C}C_g$.  To resolve the ambiguity in the sign of $\mu$, we assume that $I\in\Mp^c(V)$ is mapped to $(I,1)$.  Following \cite{rr1}, we refer to $(g,\mu)$ as the \textbf{parameters} of $a\in\Mp^c(V)$.  Note that if $a\in\Mp(V)=\ker\eta$, then the parameters $(g,\mu)$ of $a$ satisfy $\mu^2\Det_\C C_g=1$.

The unitary group $U(V)\subset\Sp(V)$ is the maximal compact subgroup of $\Sp(V)$, and consists of exactly those elements of $\Sp(V)$ that commute with the complex structure $J$.  For any $g\in U(V)$, $C_g=g$ and $\Det_\C C_g=\Det_\C g\in U(1)$.  

We view the \textbf{symplectic frame bundle} $\Sp(M,\omega)\maps{\rho}(M,\omega)$ as a right principal $\Sp(V)$ bundle over $M$, defined fiberwise such that for all $m\in M$, every $b\in\Sp(M,\omega)_m$ is a linear symplectic isomorphism $b:(V,\Omega)\rightarrow(T_mM,\omega_m)$.  The group $\Sp(V)$ acts on the fibers of $\Sp(M,\omega)$ by precomposition.

\begin{definition}\label{def:Mpc}
Let $(M,\omega)$ be a symplectic manifold with symplectic frame bundle $\Sp(M,\omega)\maps{\rho}(M,\omega)$.  A \textbf{metaplectic-c prequantization} for $(M,\omega)$ is a triple $(P,\Sigma,\gamma)$, where:
\begin{itemize}
\item $P\maps{\Pi}M$ is a right principal $\Mp^c(V)$ bundle;
\item the map $P\maps{\Sigma}\Sp(M,\omega)$ satisfies $\rho\circ\Sigma=\Pi$ and $\Sigma(q\cdot a)=\Sigma(q)\cdot\sigma(a)$ for all $q\in P$ and $a\in\Mp^c(V)$;
\item $\gamma$ is a $\mathfrak{u}(1)$-valued one-form on $P$, invariant under the principal $\Mp^c(V)$ action, such that $d\gamma=\frac{1}{i\hbar}\Pi^*\omega$, and for all $\alpha\in\mmp^c(V)$, if $\alpha$ generates the vector field $\alpha_P$ on $P$, then $\gamma(\alpha_P)=\frac{1}{2}\eta_*\alpha$.
\end{itemize}

\end{definition}

If $(P,\Sigma,\gamma)\rightarrow(M,\omega)$ is a metaplectic-c prequantization, then $(P,\gamma)\maps{\Sigma}\Sp(M,\omega)$ is a principal circle bundle with connection one-form.  The circle that acts on the fibers of $P$ is the center $U(1)\subset\Mp^c(V)$.


\section{Equivariant Metaplectic-c Prequantization}\label{sec:EquivarMpc}


\subsection{Initial constructions}\label{subsec:EquivarInit}

From now on, we assume that $(M,\omega)$ is metaplectic-c prequantizable, and we fix a metaplectic-c prequantization $(P,\Sigma,\gamma)$.  We have the bundle projection maps $P\maps{\Pi}M$ and $\Sp(M,\omega)\maps{\rho}M$.

Recall from Section \ref{subsec:HamTs} that there is a Hamiltonian $T^k$ action on $(M,\omega)$ with momentum map $\Phi$.  Each $\xi\in\mt$ generates the vector field $\xi_M$ on $M$ with flow $\phi^t_\xi$.  Since $\phi^t_\xi$ preserves $\omega$, it can be lifted to a flow $\tp_{\xi}^t$ on $\Sp(M,\omega)$, defined by $$\tp^t_{\xi}(b)=\phi^t_{\xi *}|_m\circ b,\ \ \forall m\in M,\ \forall b\in\Sp(M,\omega)_m.$$  The corresponding vector field on $\Sp(M,\omega)$ is $$\txi_M(b)=\lat{\dd{}{t}}{t=0}\tp^t_{\xi}(b).$$  Let $\tau=\exp(\xi)\in T^k$ act on $\Sp(M,\omega)$ by $\tp_\xi^1$.  It is easily verified that this definition yields a well-defined group action of $T^k$ on $\Sp(M,\omega)$ that lifts the $T^k$ action on $M$ and commutes with the principal $\Sp(V)$ action.  

Suppose there is a lift of the $T^k$ action on $\Sp(M,\omega)$ to one on $P$ that preserves $\gamma$.  Then this $T^k$ action also commutes with the principal $\Mp^c(V)$ action.  Let $\xi\in\mt$ generate the vector field $\xi_P$ on $P$.  It is immediate that $L_{\xi_P}\gamma=0$ if and only if $d(\gamma(\xi_P))=-\frac{1}{i\hbar}d\Pi^*\Phi^\xi$.  If, in particular, $$\gamma(\xi_P)=-\frac{1}{i\hbar}\Pi^*\Phi^\xi,$$ then $(P,\Sigma,\gamma)$ is called an \textbf{equivariant metaplectic-c prequantization} for $(M,\omega,\Phi)$.  An equivalent definition for spin-c structures appears in \cite{ckt}, although it is stated in terms of equivariant cohomology classes.  An analogous definition for an equivariant prequantization line bundle appears in \cite{ggk}.  

In the remainder of this section, we determine a necessary and sufficient condition for the metaplectic-c prequantized system $(P,\Sigma,\gamma)\rightarrow(M,\omega,\Phi)$ to admit an equivariant metaplectic-c prequantization.  For all $\xi\in\mt$, let $\xi_P$ be the vector field on $P$ such that $\xi_P$ is a lift of $\txi_M$ and $\gamma(\xi_P)=-\frac{1}{i\hbar}\Pi^*\Phi^\xi$.  Let $\psi^t_\xi$ be the flow of $\xi_P$ on $P$.  If there is a $T^k$ action on $P$ such that $\xi\in\mt$ generates the vector field $\xi_P$, then $\tau=\exp(\xi)\in T^k$ must act on $P$ by the map $\psi^1_\xi$.  We will find a condition that ensures that these maps $\psi^1_\xi$ yield a well-defined $T^k$ action on $P$.  It suffices to guarantee that for all $\xi\in\Z^k\subset\mt$, $\psi^1_\xi$ is the identity map on $P$.  

The following argument is based on Example 6.10 in \cite{ggk} (pp.~93-94), which establishes a necessary and sufficient condition for $(M,\omega,\Phi)$ to admit an equivariant prequantization line bundle.  Our application of their proof to a metaplectic-c prequantization requires some additional steps concerning the symplectic frame bundle.

Assume that the $T^k$ action has a fixed point $z\in M$.  For example, it is sufficient to assume that $M$ is compact.  However, noncompact examples also exist, and the remainder of the argument does not require compactness.  In Section \ref{subsec:EquivarFix}, we will determine a condition on $\Phi(z)$ such that there is a well-defined $T^k$ action on the fiber $P_z$.  Then, in Section \ref{subsec:EquivarP}, we will show that this condition guarantees a $T^k$ action on all of $P$.


\subsection{$T^k$ action over a fixed point}\label{subsec:EquivarFix}

Let $z\in M$ be a fixed point of the $T^k$ action.  Any $\tau=\exp(\xi)\in T^k$ acts on $TM$ by the pushforward $\phi^1_{\xi*}$, which preserves the symplectic form.  In particular, $\tau:T_zM\rightarrow T_zM$ is a linear symplectic isomorphism.

Let $U$ be a neighborhood of $z$ in $M$ over which $P|_U$ admits a local trivialization:  $P|_U\cong U\times\Mp^c(V)$.  This induces a local trivialization $\Sp(M,\omega)|_U\cong U\times\Sp(V)$, where the map $\Sigma$ is identified with $(\mbox{Id}_U,\sigma)$.  
$$\xymatrix{
P|_U \ar[r]^{\cong} \ar[d]^{\Sigma} & U\times\Mp^c(V) \ar[d]^{(\mbox{Id}_U,\sigma)} \\
\Sp(M,\omega)|_U \ar[r]^{\cong} & U\times\Sp(V)
}$$
It further induces a local trivialization $TM|_U\cong U\times V$, under which there is an identification of the symplectic vector space $(T_zM,\omega_z)$ with $(V,\Omega)$.  Since $T^k$ acts on $(T_zM,\omega_z)$ by linear symplectic isomorphisms, this identification yields a group homomorphism $T^k\maps{\kappa_z}\Sp(V)$.  By a suitable adjustment of the choice of trivializing section of $P|_U$, we can arrange that $$T^k\maps{\kappa_z}U(V)\subset\Sp(V).$$ For emphasis:  this property is only required to hold over the single point $z$.  On the level of tangent spaces, we obtain identifications $$T_zM\cong V,\ \ \ T_{(z,I)}\Sp(M,\omega)\cong V\times\msp(V),\ \ \ T_{(z,I)}P\cong V\times\mmp^c(V)=V\times\msp(V)\oplus\mathfrak{u}(1).$$  

Let $\xi\in\mt$ be arbitrary.  It is immediate that $\xi_M(z)=0$.  The pushforward $\phi^t_{\xi*}$, acting on $T_zM$, becomes the symplectic group element $\kappa_z(\exp(t\xi))\in\Sp(V)$.  The lifted flow $\tp^t_\xi$ on $\Sp(M,\omega)$ satisfies $$\tp^t_\xi(z,I)=(z,\kappa_z(\exp(t\xi))),$$ and so $$\txi_M(z,I)=\lat{\dd{}{t}}{t=0}\tp^t_{\xi}(z,I)=(0,\kappa_{z*}\xi).$$

The vector field $\xi_P$ is the lift of $\txi_M$ to $P$ such that $\gamma(\xi_P)=-\frac{1}{i\hbar}\Pi^*\Phi^\xi$.  At $(z,I)$, we have   $$\xi_P(z,I)=\lb{0,\kappa_{z*}\xi\oplus-\frac{1}{i\hbar}\Phi^\xi(z)}.$$ The flow $\psi^t_\xi$ of $\xi_P$ satisfies $$\psi^t_\xi(z,I)=\lb{z,\exp\ls{t\lb{\kappa_{z*}\xi\oplus-\frac{1}{i\hbar}\Phi^\xi(z)}}}.$$  The desired action of $T^k$ on $P_z$ exists if and only if $\psi^1_\xi(z,I)=(z,I)$ for all $\xi\in\Z^k\subset\mt$.

Since $U(1)$ is the center of $\Mp^c(V)$, we can write the $\Mp^c(V)$ term in the above expression for $\psi^t_\xi(z,I)$ as $$\exp\ls{t\lb{\kappa_{z*}\xi\oplus-\frac{1}{i\hbar}\Phi^\xi(z)}}=\exp\lb{t\kappa_{z*}\xi\oplus 0}e^{2\pi it\Phi^\xi(z)/h},$$ where $h=2\pi\hbar$, and where $\exp\lb{t\kappa_{z*}\xi\oplus 0}\in\Mp(V)\subset\Mp^c(V)$ and $e^{2\pi it\Phi^\xi(z)/h}\in U(1)\subset\Mp^c(V)$.

The parameters of $\exp(t\kappa_{z*}\xi\oplus 0)\in\Mp(V)$ take the form $(\kappa_z(\exp(t\xi)),\mu(t))$ where $$\mu(t)^2\Det_\C C_{\kappa_z(\exp(t\xi))}=1.$$  Note that $t=0$ corresponds to $I\in\Mp(V)$, so we must have $\mu(0)=1$.  Further, since $\kappa_z(\exp(t\xi))\in U(V)$, we have $C_{\kappa_z(\exp(t\xi))}=\kappa_z(\exp(t\xi))$ and $\Det_\C \kappa_z(\exp(t\xi))\in U(1)$.  Let $\Delta$ denote the map from $\Sp(V)$ to $\C$ given by $g\mapsto\Det_\C C_g$.  Let $\Delta$ also denote the restriction of this map to $U(V)$, where it is given by $g\mapsto\Det_\C g$, and note that $U(V)\maps{\Delta}U(1)$ is a group homomorphism. 
$$\xymatrix{
T^k \ar[r]^{\kappa_z} \ar[rd]^{w_z} & U(V) \ar[r]^{\subset} \ar[d]^{\Delta=\Det_\C} & \Sp(V) \ar[d]^{\Delta=\Det_\C C_g} \\
 & U(1) \ar[r]^{\subset} & \C
}$$
Let $w_z=\Delta\circ\kappa_z$.  Then $$\mu(t)^2w_z(\exp(t\xi))=1$$ and $\mu(0)=0$, which implies that $\mu(t)=w_z(\exp(t\xi))^{-1/2}$.  Thus the parameters of $\exp(t\kappa_{z*}\xi\oplus 0)$ are $$(\kappa_z(\exp(t\xi)),w_z(\exp(t\xi))^{-1/2}).$$

Identify $\mathfrak{u}(1)$ with $\R$ such that for all $\lambda\in\mathfrak{u}(1)$, $\exp(\lambda)=e^{2\pi i\lambda}\in U(1)$.  Then $w_z(\exp(t\xi))=e^{2\pi it w_{z*}\xi}\in U(1)$, and so $w_z(\exp(t\xi))^{-1/2}=e^{-\pi it w_{z*}\xi}$.  Thus the parameters of $\exp(t\kappa_{z*}\xi\oplus 0)$ are $$(\kappa_z(\exp(t\xi)),e^{-\pi itw_{z*}\xi}),$$ which implies that the parameters of $\exp\ls{t\lb{\kappa_{z*}\xi\oplus-\frac{1}{i\hbar}\Phi^\xi(z)}}$ are $$\lb{\kappa_z(\exp(t\xi)),e^{-\pi itw_{z*}\xi}e^{2\pi it\Phi^\xi(z)/h}}.$$  

Now assume that $\xi\in\Z^k\subset\mt$, and set $t=1$.  The condition $\psi^1_\xi(z,I)=(z,I)$ is satisfied if and only if $$\lb{\kappa_z(\exp(\xi)),e^{-\pi iw_{z*}\xi}e^{2\pi i\Phi^\xi(z)/h}}=(I,1).$$ 
It is clear that $\kappa_z(\exp(\xi))=I$, and it remains to ensure that $e^{-\pi iw_{z*}\xi}e^{2\pi i\Phi^\xi(z)/h}=1$.  This equation holds if and only if $$-\pi iw_{z*}\xi+\frac{2\pi i\Phi^\xi(z)}{h}=2\pi i N$$ for some $N\in\Z$, which rearranges to $$\frac{1}{h}\Phi^\xi(z)-\frac{1}{2}w_{z*}\xi=N.$$  Since an equation of this form must hold for all $\xi\in\Z^k\subset\mt$, we conclude that the value $\Phi(z)\in\mt^*$ must satisfy $$\frac{1}{h}\Phi(z)-\frac{1}{2}w_{z*}\in\Z^{k*}\subset\mt^*.$$  This is similar to, but shifted from, the result in \cite{ggk} that an equivariant prequantization line bundle exists if and only if $\frac{1}{h}\Phi(z)$ is in the integer lattice $\Z^{k*}\subset\mt^*$ (adjusted for differing $\hbar$ conventions).


\subsection{$T^k$ action on $P$}\label{subsec:EquivarP}

We continue to follow a modified version of the argument in Example 6.10 of \cite{ggk}.  Fix $\xi\in\Z^k\subset\mt$.  Then $\psi^1_\xi$ is a lift of the identity maps $\tp^1_\xi$ on $\Sp(M,\omega)$ and $\phi^1_\xi$ on $M$.  
\begin{itemize}
\item Since $\psi^1_\xi$ is a lift of the identity map on $M$, there is a map $R_\xi:M\rightarrow\Mp^c(V)$ such that $\psi^1_\xi(q)=q\cdot R_\xi(\Pi(q))$ $\forall q\in P.$
\item Since $\psi^1_\xi$ is a lift of the identity map on $\Sp(M,\omega)$, there is a map  $\widetilde{R}_\xi:\Sp(M,\omega)\rightarrow U(1)\subset\Mp^c(V)$ such $\psi^1_\xi(q)=q\cdot\widetilde{R}_\xi(\Sigma(q))$ $\forall q\in P$.
\end{itemize}
These observations together imply that the target of the map $R_\xi$ is $U(1)$.  That is, there is a map $R_\xi:M\rightarrow U(1)$ such that $\psi^1_\xi(q)=q\cdot R_\xi(\Pi(q))$ for all $q\in P$.

Assume that the condition on $\Phi$ derived in the previous section has been satisfied over the fixed point $z$.  Then $R_\xi(z)=1$.  It remains to show that $R_\xi$ is constant over $M$.  Let $u(s)$ be an arbitrary path in $M$, where $s\in[0,1]$.  We will show that $R_\xi$ is constant over $u(s)$.

Recall that the determinant map $\Mp^c(V)\maps{\eta}U(1)$ acts on $\lambda\in U(1)\subset\Mp^c(V)$ by $\eta(\lambda)=\lambda^2$.  Let $Y\maps{\pi}M$ be the circle bundle associated to $P\maps{\Pi}M$ by $\eta$, and let $P\maps{H}Y$ be the corresponding bundle map.  It is easily verified that $\gamma$ is basic with respect to $H$.  Let $\gamma^\eta$ be the $\mathfrak{u}(1)$-valued one-form on $Y$ such that $H^*\gamma^\eta=2\gamma$.  Then $\gamma^\eta$ is a connection one-form on $Y$, and $d\gamma^\eta=\frac{2}{i\hbar}\pi^*\omega$.

The pushforward $H_*\xi_P$ is a well-defined vector field on $Y$, which we denote by $\xi_Y$.   Let the flow of $\xi_Y$ on $Y$ be $\hchi^t_\xi$.  The various vector fields and their flows are summarized below.
$$\xymatrix{
\xi_P,\psi^t_\xi & (P,\gamma) \ar[r]^H \ar[d] & (Y,\gamma) \ar[dd] & \xi_Y,\hchi^t_\xi \\
\txi_M,\tp^t_\xi & \Sp(M,\omega) \ar[d] & & \\
\xi_M,\phi^t_\xi & (M,\omega) \ar[r]^{=} & (M,\omega) & \xi_M,\phi^t_\xi.
}$$
Since $H_*\xi_P=\xi_Y$, $H$ intertwines the flows $\psi^t_\xi$ and $\hchi^t_\xi$.  In particular, $$\hchi^1_\xi\circ H(q)=H\circ\psi^1_\xi(q)=H(q\cdot R_\xi(\Pi(q)))=H(q)\cdot R^2_\xi(\Pi(q)).$$  That is, the map $\hchi^1_\xi$ acts on $Y$ by $$\hchi^1_\xi(y)=y\cdot R^2_\xi(\pi(y)),\ \ \forall y\in Y.$$

From the definitions of $\xi_Y$ and $\gamma^\eta$, it follows that $\gamma^\eta(\xi_Y)=-\frac{2}{i\hbar}\pi^*\Phi^\xi$.  Let $\partial_\theta$ be the vertical vector field on $Y$ such that $\gamma^\eta(\partial_\theta)=2\pi i\in\mathfrak{u}(1)$.  Then we have $$\xi_Y=(\xi_M)_{hor}+\frac{2}{h}\Phi^\xi\partial_\theta,$$ where $(\xi_M)_{hor}$ represents the lift of $\xi_M$ to $Y$ that is horizontal with respect to $\gamma^\eta$.

Let $C=\R/\Z\times[0,1]$, with coordinates $(r,s)$.  Define $F:C\rightarrow M$ by $$F(r,s)=\exp(r\xi)\cdot u(s),\ \ \forall(r,s)\in C,$$ and let $\omega_C=F^*\omega$.  Construct the pullback of $(Y,\gamma^\eta)$ to $C$:
$$\xymatrix{
(D,\delta) \ar[r]^{\hF} \ar[d] & (Y,\gamma^\eta) \ar[d]^{\pi} \\
(C,\omega_C) \ar[r]^F & (M,\omega)
}$$
where $\delta=\hF^*\gamma^\eta$ and $$D=F^*Y=\cb{(c,y)\in C\times Y:F(c)=\pi(y)}.$$  The bundle map $\hF:D\rightarrow Y$ acts by $$\hF(c,y)=y,\ \ \forall (c,y)\in D.$$   By construction, $(D,\delta)$ is a circle bundle with connection one-form over $(C,\omega_C)$.  Let $\partial_\theta$ also denote the vertical vector field on $D$ such that $\delta(\partial_\theta)=2\pi i\in\mathfrak{u}(1)$, and note that $\hF_*\partial_\theta=\partial_\theta$.

Abbreviate the vector fields $\pp{}{r}$ and $\pp{}{s}$ on $C$ by $\partial_r$ and $\partial_s$, with flows $\psi^t_r$ and $\psi^t_s$ respectively.  It is immediate from the definition of $F$ that $$F_*|_c\partial_r=\xi_M(F(c)),\ \ \forall c\in C.$$

Let $\Psi:C\rightarrow\R$ be given by $$\Psi(r,s)=\Phi^\xi(u(s)),\ \ \forall(r,s)\in C.$$  We claim that $\omega_C=\pp{\Psi}{s}dr\wedge ds$.  This is established by calculating, at arbitrary $c=(r,s)\in C$, $$\partial_r\contr(\omega_C)_c=F^*d\Phi^\xi_{F(c)},$$ and $$\partial_s\contr F^*d\Phi^\xi_{F(c)}=\lat{\pp{\Psi}{s}}{(r,s)}.$$

Let $\zeta_r$ be the vector field on $D$ given by $$\zeta_r=(\partial_r)_{hor}+\frac{2}{h}\Psi\partial_\theta,$$ with flow $\hpsi^t_r$.  Then for all $(c=(r,s),y)\in D$, $$\hF_*|_{(c,y)}\zeta_r=(\xi_M)_{hor}(\hF(c,y))+\frac{2}{h}\Phi^\xi(u(s))\partial_\theta=\xi_Y(\hF(c,y)).$$  Thus $\hF$ intertwines the flows $\hpsi^t_r$ of $\zeta_r$ and $\hchi^t_\xi$ of $\xi_Y$.  In particular, at $t=1$, $$\hF\circ\hpsi^1_r(c,y)=\hchi^1_\xi\circ\hF(c,y)=\hchi^1_\xi(y)=y\cdot R^2_\xi(\pi(y))=\hF(c,y\cdot R^2_\xi(\pi(y))),$$ which implies that $$\hpsi^1_r(c,y)=(c,y\cdot R^2_\xi(\pi(y)))=(c,y)\cdot R^2_\xi(\pi(y)).$$  That is, $\hpsi^1_r$ fixes the base $C$ and rotates each fiber $D_{(c,y)}$ by $R^2_\xi(\pi(y))$.  Note that $\pi(y)=F(c)$, so we have $$R^2_\xi(\pi(y))=R^2_\xi(F(c))=R^2_\xi(\exp(r\xi)u(s)).$$

Let $\zeta_s=(\partial_s)_{hor}$ on $D$, with flow $\hpsi^t_s$.  A standard calculation establishes that $[\zeta_r,\zeta_s]=0$.  Hence their flows $\hpsi^t_r$ and $\hpsi^t_s$ commute.  In particular, for any $(c=(r,s),y)\in D$,  $$\hpsi^1_r\circ\hpsi^t_s(c,y)=\hpsi^t_s\circ\hpsi^1_r(c,y)=\hpsi^t_s(c,y\cdot R^2_\xi(\pi(y)))=\hpsi^t_s(c,y)\cdot R^2_\xi(\exp(r\xi)\cdot u(s)),$$ having used the fact that $\hpsi^t_s$ is the flow of a horizontal vector field on $D$ and is therefore equivariant with respect to the principal $U(1)$ action.  We also calculate 
\begin{align*}
\hpsi^1_r\circ\hpsi^t_s(c,y)=&\hpsi^t_s(c,y)\cdot R^2_\xi(\pi(\hF(\hpsi^t_s(c,y))))=\hpsi^t_s(c,y)\cdot R^2_\xi(F(\psi^t_s(c)))\\
=&\hpsi^t_s(c,y)\cdot R^2_\xi(\exp(r\xi)\cdot u(s+t)),
\end{align*}
where we note that $\psi^t_s(c)=(r,s+t)$.  We conclude that $R^2_\xi(\exp(r\xi)\cdot u(s))=R^2_\xi(\exp(r\xi)\cdot u(s+t))$ for all $r,s,t$.  Hence $R^2_\xi$ is constant over the path $u$.  Since $u$ was arbitrary, $R^2_\xi$ is in fact constant on $M$.

Recall that $R_\xi(z)=1$.  Thus $R^2_\xi=1$ everywhere on $M$, and consequently $R_\xi=1$ everywhere on $M$, as needed.  Hence the $T^k$ action is well defined everywhere on $P$.  We have now proved the following theorem.

\begin{theorem}\label{thm:Equivar}
Let $(M,\omega)$ have an effective Hamiltonian $T^k$ action with momentum map $\Phi$ and a fixed point $z$.  Then $(M,\omega,\Phi)$ admits an equivariant metaplectic-c prequantization if and only if $(M,\omega)$ is metaplectic-c prequantizable and the momentum map $\Phi$ satisfies $$\frac{1}{h}\Phi(z)-\frac{1}{2}w_{z*}\in\Z^{k*}\subset\mt^*.$$
\end{theorem}

Assume the hypotheses of the theorem, and let $z$ be a fixed point for the $T^k$ action on $(M,\omega)$.  By adding a constant to the momentum map $\Phi$, it is always possible to satisfy the condition $\frac{1}{h}\Phi(z)-\frac{1}{2}w_{z*}\in\Z^{k^*}$.  Thus any metaplectic-c prequantization $(P,\Sigma,\gamma)$ for $(M,\omega)$ can be converted to an equivariant metaplectic-c prequantization by a suitable shift of the momentum map.


\section{Quantized Energy Levels}\label{sec:QuantE}

In this section, we extend the concept of a quantized energy level to the equivariant metaplectic-c prequantized system $(P,\Sigma,\gamma)\rightarrow(M,\omega,\Phi)$.  Section \ref{subsec:PS} reviews the constructions due to Robinson \cite{rob1} that we use to define a quantized energy level, and concludes with the generalization of our definition from \cite{v2} to a regular value of the momentum map $\Phi$.  In Section \ref{subsec:GenQuantE}, we determine the quantized energy levels of the system $(M,\omega,\Phi)$, assuming that $\Phi$ has been shifted so that the metaplectic-c prequantization is equivariant.


\subsection{Descending to the quotient}\label{subsec:PS}

Assume that $(M,\omega)$ is metaplectic-c prequantizable, and let $(P,\Sigma,\gamma)$ be a metaplectic-c prequantization: $$(P,\gamma)\maps{\Sigma}\Sp(M,\omega)\maps{\rho}(M,\omega),\ \ \rho\circ\Sigma=\Pi,\ \ d\gamma=\frac{1}{i\hbar}\Pi^*\omega.$$  Assume that there is an effective Hamiltonian $T^k$ action on $M$ with momentum map $\Phi$ and at least one fixed point.  Shift $\Phi$ if necessary so that $(P,\Sigma,\gamma)$ is an equivariant metaplectic-c prequantization for $(M,\omega,\Phi)$.

Let $x\in\mt^*$ be a regular value of the momentum map $\Phi$, and let $S=\Phi^{-1}(x)$.  Then $S$ is a codimension-$k$ embedded submanifold of $M$.  Recall that $\cb{\xi_1,\ldots,\xi_k}$ is the standard basis for $\R^k$.  For all $s\in S$, the symplectic orthogonal to $T_sS$ is $T_s S^\perp=\mbox{span}\cb{\xi_{1M}(s),\ldots,\xi_{kM}(s)}\subset T_sS$, implying that $S$ is a coisotropic submanifold of $M$.  In the special case where $k=n$, $S$ is a Lagrangian submanifold.

Within the model symplectic vector space $(V,\Omega)$, let $W\subset V$ be a coisotropic subspace of codimension $k$, with symplectic orthogonal $W^\perp\subset W$.  Then $W/W^\perp$ is a symplectic vector space with a symplectic structure inherited from $V$.  Let $\Sp(V;W)\subset\Sp(V)$ be the subgroup that preserves $W$.  There is a natural group homomorphism $\Sp(V;W)\maps{\nu}\Sp(W/W^\perp)$.  

Let $\Mp^c(V;W)\subset\Mp^c(V)$ be the preimage of $\Sp(V;W)$ under $\sigma$, and let $\Mp^c(W/W^\perp)$ be the metaplectic-c group for $W/W^\perp$.  Robinson and Rawnsley \cite{rr1} showed that $\nu$ lifts to a group homomorphism $\hnu$ on the level of metaplectic-c groups:
$$\xymatrix{
\Mp^c(V) \ar[d]^{\sigma} & \Mp^c(V;W) \ar[l]_{\supset} \ar[r]^{\hnu} \ar[d]^{\sigma} & \Mp^c(W/W^\perp) \ar[d]^{\sigma} \\
\Sp(V) & \Sp(V;W) \ar[l]_{\supset} \ar[r]^{\nu} & \Sp(W/W^\perp)
}$$
The lifted map $\hnu$ has the property that $\hnu_*\circ\eta_*=\eta_*$.

In the following construction, which is due to Robinson \cite{rob1}, the above diagram serves as a model for fiberwise constructions over $(M,\omega)$.  The first column corresponds to the original three-level structure, $$(P,\gamma)\rightarrow\Sp(M,\omega)\rightarrow(M,\omega).$$  For the second column, let $\Sp(M,\omega;S)\subset\Sp(M,\omega)$ be the subbundle, lying only over $S$, such that for all $s\in S$ and all $b\in\Sp(M,\omega;S)_s$, $bW=T_sS$.  Lastly, let $(P^S,\gamma^S)$ be the pullback of $(P,\gamma)$ to $\Sp(M,\omega;S)$.  

For the third column, treat $W/W^\perp$ as a model symplectic vector space for the symplectic vector bundle $TS/TS^\perp\rightarrow S$, so that the symplectic frame bundle $\Sp(TS/TS^\perp)\rightarrow S$ becomes a right principal $\Sp(W/W^\perp)$ bundle over $S$.  Then $\Sp(TS/TS^\perp)$ is naturally identified with the bundle associated to $\Sp(M,\omega;S)$ by the map $\nu$.  Let $P_S$ be the bundle associated to $P^S$ by the map $\hnu$.  The properties of $\hnu$ guarantee that there is a one-form $\gamma_S$ on $P_S$ such that $\hnu^*\gamma_S=\gamma^S$.  Then $(P_S,\gamma_S)\rightarrow\Sp(TS/TS^\perp)$ is a principal circle bundle with connection one-form.

$$\xymatrix{
(P,\gamma) \ar[d] & (P^S,\gamma^S) \ar[d] \ar[l]_{\supset} \ar[r]^{\hnu} & (P_S,\gamma_S) \ar[d] \\
\Sp(M,\omega) \ar[d] & \Sp(M,\omega;S) \ar[d] \ar[l]_{\supset} \ar[r]^{\nu} & \Sp(TS/TS^\perp) \ar[d] \\
(M,\omega) & S \ar[l]_{\supset} \ar[r]^{=} & S
}$$

Let $\xi\in\mt$ be arbitrary.  Construct $\xi_M$ with flow $\phi^t_\xi$ on $M$, $\txi_M$ with flow $\tp^t_\xi$ on $\Sp(M,\omega)$, and $\xi_P$ with flow $\psi^t_\xi$ on $P$, as described in Section \ref{subsec:EquivarInit}.  Recall in particular that $\xi_P$ is the lift of $\txi_M$ to $P$ such that $\gamma(\xi_P)=-\frac{1}{i\hbar}\Pi^*\Phi^\xi$.  Let $\hxi_M$ be the lift of $\txi_M$ to $P$ that is horizontal with respect to $\gamma$, and let its flow be $\hp^t_\xi$.  Each of the vector fields $\xi_M$, $\txi_M$, $\hxi_M$ and $\xi_P$ restricts to a vector field on the appropriate manifold in the second column, and descends to a vector field on the manifold in the third column.  Moreover, since the $T^k$ action commutes with all of the principal actions, the $T^k$ action on each level preserves the manifold in the second column, and descends to a $T^k$ action on the manifold in the third column.

Viewing $(P,\gamma)$ as a circle bundle with connection one-form over $\Sp(M,\omega)$, let $\partial_\theta$ be the vertical vector field on $P$ such that $\gamma(\partial_\theta)=2\pi i\in\mathfrak{u}(1)$.  Then $\xi_P=\hxi_M+\frac{1}{h}\Phi^\xi\partial_\theta$.  We also denote by $\partial_\theta$ the restriction of this vector field to $P^S$, and the induced vertical vector field on $P_S$.  In each column, we have $\gamma(\partial_\theta)=\gamma^S(\partial_\theta)=\gamma_S(\partial_\theta)=2\pi i\in\mathfrak{u}(1)$.  Note that for all $s\in S$, $\Phi^\xi(s)=x\cdot\xi$, since $S$ is the $x$-level set of $\Phi$.  Thus, in columns 2 and 3, $$\xi_P=\hxi_M+\frac{1}{h}(x\cdot\xi)\partial_\theta.$$


\subsection{Generalized quantized energy levels}\label{subsec:GenQuantE}

In our previous paper \cite{v2}, we considered the system $(M,\omega,H)$, where $H\in C^\infty(M)$ is viewed as a Hamiltonian energy function.  Let $(P,\Sigma,\gamma)\rightarrow(M,\omega)$ be a metaplectic-c prequantization.  Let $E$ be a regular value of $H$, and use the embedded coisotropic submanifold $S=H^{-1}(E)\subset M$ to construct the three columns of three-level structures as described in the previous section.  Let $\xi_H$ be the Hamiltonian vector field for $H$ on $M$, and let its lift to $\Sp(M,\omega)$ be $\txi_H$.  Then $\xi_H$ and $\txi_H$ restrict to column 2 and descend to column 3.  We define $E$ to be a \textbf{quantized energy level} for the system $(M,\omega,H)$ if $\gamma_S$ has trivial holonomy over all closed orbits of $\txi_H$ on $\Sp(TS/TS^\perp)$.  

This definition has a natural generalization to a regular value of a family of $k$ Poisson-commuting functions.  Recall that $\cb{\xi_1,\ldots,\xi_k}$ is the standard basis for $\mt=\R^k$.  In our context, the family of Poisson-commuting functions consists of the $k$ components of the momentum map, $\Phi^{\xi_1},\ldots,\Phi^{\xi_k}$.  The generalized definition is as follows.

\begin{definition}\label{def:QuantE}
Let $(P,\Sigma,\gamma)$ be an equivariant metaplectic-c prequantization for $(M,\omega,\Phi)$.  Let $x$ be a regular value of $\Phi$, and let $S=\Phi^{-1}(x)$.  Let $F$ be the distribution on $\Sp(TS/TS^\perp)$ spanned by the vector fields $\left\{\txi_{1M},\right.\left.\ldots,\txi_{kM}\right\}.$  Then $x$ is a \textbf{quantized energy level} for $(M,\omega,\Phi)$ if $\gamma_S$ has trivial holonomy over all of the leaves of $F$.
\end{definition}

The connection one-form $\gamma_S$ is flat over each leaf of the distribution spanned by $\cb{\txi_{1M},\ldots,\txi_{kM}}$.  For the regular value $x$ to be a quantized energy level, it suffices to ensure that $\gamma_S$ has trivial holonomy over the orbits of each $\txi_{jM}$.  Given any initial point $b\in\Sp(TS/TS^\perp)$, the integral curve $\tp^t_{\xi_j}(b)$ satisfies $\tp^1_{\xi_j}(b)=b$.  We need to show that the integral curves $\hp^t_{\xi_j}(q)$ of the horizontal lift $\hxi_{jM}$ on $P_S$ satisfy $\hp^1_{\xi_j}(q)=q$, for all $q\in P_S$.

As previously noted, the vector fields $\xi_{jP}$ and $\hxi_{jM}$ on $P_S$ are related by $$\xi_{jP}=\hxi_{jM}+\frac{1}{h}(x\cdot\xi_j)\partial_\theta.$$  Since these two vector fields differ by a constant multiple of $\partial_\theta$ everywhere on $P_S$, their flows are related by $$\psi^t_{\xi_j}(q)=\hp^t_{\xi_j}(q)\cdot e^{2\pi it(x\cdot\xi_j)/h},\ \ \forall q\in P.$$  

We know that $\tau=\exp(\xi)\in T^k$ acts on $P_S$ by the map $\psi^1_{\xi}$ for all $\xi\in\mt$, and this is a well-defined $T^k$ action.  In particular, $\exp(\xi_j)=I\in T^k$, so $\psi^1_{\xi_j}$ is the identity map on $P_S$.  We have $$\psi^1_{\xi_j}(q)=q=\hp^1_{\xi_j}(q)e^{2\pi i(x\cdot\xi_j)/h}.$$  Thus $\hp^1_{\xi_j}(q)=q$ if and only if $e^{2\pi i(x\cdot\xi_j)/h}=1$, which occurs when $x\cdot\xi_j=N_jh$ for some $N_j\in\Z$.  This condition is satisfied for all $\xi_j$ exactly when $x\in h\Z^{k^*}\subset\mt^*$.  We conclude the following result.

\begin{theorem}\label{thm:EquivarQuantE}
Let $(M,\omega)$ have an effective Hamiltonian $T^k$ action with momentum map $\Phi$ and at least one fixed point.  Assume that $(M,\omega)$ is metaplectic-c prequantizable, and shift $\Phi$ by a constant if necessary so that $(M,\omega,\Phi)$ admits an equivariant metaplectic-c prequantization.  Then the regular value $x\in\mt^*$ for $\Phi$ is a quantized energy level for the system if and only if $$x\in h\Z^{k*}\subset\mt^*.$$
\end{theorem}

An immediate consequence of this theorem arises in the context of symplectic reduction.  Assume the hypotheses of Theorem \ref{thm:EquivarQuantE}, and let $x\in\mt^*$ be a regular value of $\Phi$.  Further assume that $T^k$ acts freely on the level set $S=\Phi^{-1}(x)$.  Then, by Marsden-Weinstein reduction, the space of orbits $M_0=S/T^k$ is a manifold, and it acquires a symplectic form $\omega_0$.  Let $S\maps{\varrho}M_0$ be the quotient map from $S$ to its orbit space.  Then $\varrho^*\omega_0=\omega_S$, where $\omega_S$ is the pullback of $\omega$ to $S$.

Let $s\in S$ be arbitrary and let $\varrho(s)=m_0$.  On the level of tangent spaces, we have the short exact sequence $$0\rightarrow T_sS^\perp\rightarrow T_sS\maps{\varrho_*|_s} T_{m_0}M_0\rightarrow 0,$$ implying that $\varrho_*|_s$ induces a linear symplectic isomorphism between $T_sS/T_sS^\perp$ and $T_{m_0}M_0$.  For all $\xi\in\mt$, the pushforward $\phi^t_{\xi*}:TM\rightarrow TM$ preserves both $TS$ and $TS^\perp$, and so descends to a map on $TS/TS^\perp$.  If we let $\tau=\exp(\xi)\in T^k$ act on $TS/TS^\perp$ by the pushforward $\phi^1_{\xi*}$, the result is a $T^k$ action on $TS/TS^\perp$.  Using these observations, it is straightforward to verify that the tangent bundle $TM_0$ is naturally identified with the quotient $(TS/TS^\perp)/T^k$.  From this it follows that the symplectic frame bundle $\Sp(M_0,\omega_0)$ is naturally identified with the quotient $\Sp(TS/TS^\perp)/T^k$.

$$\xymatrix{
(P_S,\gamma_S) \ar[d] & \\
\Sp(TS/TS^\perp) \ar[d] \ar[r]_{/T^k} & \Sp(M_0,\omega_0) \ar[d] \\
S \ar[r]^{\varrho}_{/T^k} & (M_0,\omega_0)
}$$

The following fact was stated by Robinson \cite{rob1}.  Let $(Y,\gamma)\rightarrow Z$ be a principal circle bundle with connection one-form over an arbitrary manifold $Z$.  Let $F$ be a fibrating foliation of $Z$ with leaf space $Z_0$ and quotient map $Z\maps{\varrho}Z_0$.  Let the curvature of $\gamma$ be $\varpi$.  If $\gamma$ has trivial holonomy over the leaves of $F$, then $(Y,\gamma)$ descends to a principal circle bundle with connection one-form $(Y_0,\gamma_0)\rightarrow Z_0$ such that the curvature $\varpi_0$ of $\gamma_0$ satisfies $\varrho^*\varpi_0=\varpi$.

In our case, the base manifold is $\Sp(TS/TS^\perp)$, and the fibrating foliation is $F=\mbox{span}\left\{\txi_{1M},\ldots,\right.$ $\left.\txi_{kM}\right\}$.  By Definition \ref{def:QuantE}, if $x\in\mt^*$ is a quantized energy level for $(M,\omega,\Phi)$, then $\gamma_S$ has trivial holonomy over the leaves of $F$.  This is exactly the condition required for $(P_S,\gamma_S)\rightarrow\Sp(TS/TS^\perp)$ to descend to a circle bundle $$(P_0,\gamma_0)\maps{\Sigma_0}\Sp(M_0,\omega_0),$$ where $P_0=P_S/T^k$ and $\gamma_0$ is a connection one-form on $P_0$ such that the curvature of $\gamma_S$ is the pullback of the curvature of $\gamma_0$.  It is easily checked that $(P_0,\Sigma_0,\gamma_0)$ is a metaplectic-c prequantization for $(M_0,\omega_0)$.  In other words, when $x$ is a quantized energy level for $(M,\omega,\Phi)$, the top row of the diagram above can be completed in the obvious manner, and the result is a metaplectic-c prequantization for the symplectic reduction.

By Theorem \ref{thm:EquivarQuantE}, the quantized energy levels of $(M,\omega,\Phi)$ are the regular values of $\Phi$ that lie in $h\Z^{k*}$.  We conclude the following.

\begin{theorem}\label{thm:MpcRed}
Let $(M,\omega)$ have an effective Hamiltonian $T^k$ action with momentum map $\Phi$ and at least one fixed point.  Assume that $(M,\omega)$ is metaplectic-c prequantizable, and let $(P,\Sigma,\gamma)$ be an equivariant metaplectic-c prequantization for $(M,\omega,\Phi)$.  If $x\in\mt^*$ is a regular value of $\Phi$ that lies in $h\Z^{k*}$, and if $T^k$ acts freely on the level set $S=\Phi^{-1}(x)$, then the symplectic reduction $(M_0,\omega_0)$ of this level set acquires a metaplectic-c prequantization by taking the quotient of $(P_S,\gamma_S)$ by $T^k$.
\end{theorem}


\section{Examples}\label{sec:Ex}


\subsection{Harmonic oscillators}\label{subsec:SHO}

Let $M=\R^{2n}=\C^n$, with Cartesian coordinates $(q_1,\ldots,q_n,p_1,\ldots,p_n)$ and complex coordinates $z_j=q_j+ip_j$, $1\leq j\leq n$.  Equip $M$ with the symplectic form $\omega=\sum_{j=1}^ndq_j\wedge dp_j=\frac{1}{2i}\sum_{j=1}^n d\oz_j\wedge dz_j$.  Since $M$ is contractible, $(M,\omega)$ admits a metaplectic-c prequantization that is unique up to isomorphism.

Let the circle $T^1=U(1)$ act on $M$ as follows:  given $\tau\in T^1$ and $m=(z_1,\ldots,z_n)\in M$, $$\tau\cdot m=(\tau z_1,\ldots,\tau z_n).$$  Identify $\mt$ with $\R$ such that for all $\xi\in\mt$, $\exp(\xi)=e^{2\pi i\xi}\in T^1$.  The momentum map $\Phi:M\rightarrow\mt^*$ is given, up to an additive constant, by $$\Phi(m)=-\pi\sum_{j=1}^n|z_j|^2,\ \ \forall m\in M.$$

The fixed point of the $T^1$ action is the origin, $z=(0,\ldots,0)$.  If we identify $T_zM$ with $\R^{2n}=\C^n$ in the obvious way, it is immediate that every $\tau\in T^1$ acts on $T_zM$ as a complex linear isomorphism.  Explicitly, for any $\xi\in\mt$, if $\tau=\exp(\xi)$, then $\tau$ acts on $T_zM$ by the complex matrix $\mbox{diag}\lb{e^{2\pi i\xi},\ldots,e^{2\pi i\xi}}$.  The group homomorphism $w_z:T^k\rightarrow U(1)$, as defined in Section \ref{subsec:EquivarFix}, is $$w_z(\tau)=\Det_\C \mbox{diag}\lb{e^{2\pi i\xi},\ldots,e^{2\pi i\xi}}=e^{2\pi in\xi}.$$  Therefore $$\frac{1}{2}w_{z*}=\frac{n}{2}\in\mt^*.$$

Over the fixed point $z$, we find that $$\frac{1}{h}\Phi(z)-\frac{1}{2}w_{z*}=-\frac{1}{2}w_{z*}=-\frac{n}{2}.$$ If $n$ is even, then $-\frac{n}{2}\in\Z^{*}\subset\mt^*$ and the equivariance condition is satisfied, but not if $n$ is odd.  Let $$\Phi'=\Phi+\frac{hn}{2}.$$  Then $\Phi'$ is also a momentum map for the $T^1$ action, and $\frac{1}{h}\Phi'(z)-\frac{1}{2}w_*\in\Z^*\subset\mt^*$, for all $n$.

The quantized energy levels of the system $(M,\omega,\Phi')$ are the regular values of $\Phi'$ that are in $h\Z^{*}\subset\mt^*$.  Since $$\Phi'(m)=-\pi\sum_{j=1}^n|z_j|^2+\frac{hn}{2},$$ a regular value is $x\in\R$ such that $x<\frac{hn}{2}$.  Thus a quantized energy level takes the form $$x=-hN,$$ where $N\in\Z$ is such that $N>-\frac{n}{2}$.

The Hamiltonian energy function for an $n$-dimensional harmonic oscillator is $$H=\frac{1}{2}\sum_{j=1}^n(q_j^2+p_j^2)=\frac{1}{2}\sum_{j=1}^n|z_j|^2.$$  Note that $H=-\frac{1}{2\pi}\Phi'+\frac{\hbar n}{2}$.  Therefore the quantized energy levels for the system $(M,\omega,H)$ take the form $$E=\hbar\lb{N+\frac{n}{2}},$$ where $N\in\Z$ is such that $N>-\frac{n}{2}$.  

For comparison, the standard quantum mechanical calculation for the energy levels of the quantized harmonic oscillator yields $$E=\hbar\lb{N+\frac{n}{2}},$$ where $N\in\Z$ is such that $N\geq 0$.  The two calculations do not agree on the starting point for the energy spectrum (but see below).  However, the equivariant metaplectic-c prequantization does yield the $\frac{n}{2}$ shift in the energy levels.  By contrast, Kostant-Souriau quantization requires the half-form correction to obtain this shift, which uses a choice of polarization.

In quantum mechanics, the quantized energy levels of this system are obtained by solving the Schr\"odinger equation, which is linear:  an $n$-dimensional harmonic oscillator is equivalent to $n$ independent one-dimensional harmonic oscillators.  Consider the system described by the functions $H_1,\ldots,H_n$ where $H_j=\frac{1}{2}(q_j^2+p_j^2)$ for each $j$.  By an essentially identical calculation, we find that the quantized energy levels for such a system have the form $(E_1,\ldots,E_n)$, where for each $j$, $E_j=\hbar\lb{N_j+\frac{1}{2}}$ for some $N_j\in\Z$ such that $N_j>-\frac{1}{2}$.  If we view the quantized energy levels of the $n$-dimensional harmonic oscillator as $E=E_1+\ldots+E_n$, we obtain$$E=\hbar\lb{N+\frac{n}{2}},\ \ N\in\Z,\ \ N\geq 0.$$  We now have both the $\frac{n}{2}$ shift and the correct starting point, suggesting that this is the more appropriate mathematical interpretation of the physical system.


\subsection{Complex projective space}\label{subsec:CP}

Consider $\C^{n+1}$ with the usual complex coordinates $z=(z_0,\ldots,z_n)$, and complex projective space $\CP^n$ with the usual homogeneous coordinates $[z]=[z_0:\ldots:z_n]$.  The two-form $$\varpi_{FS}=\partial\overline{\partial}\log\lb{|z|^2}$$ on $\C^{n+1}$ descends to $\CP^n$.  Let $$\omega_{FS}=Ki\varpi_{FS}$$ on $\CP^n$, where $K>0$ is a positive constant to be determined.  Then $\omega_{FS}$ is a K\"ahler form on $\CP^n$:  specifically, a scalar multiple of the Fubini-Study form.  

Robinson and Rawnsley \cite{rr1} demonstrated that $(\CP^n,\omega_{FS})$ admits a metaplectic-c prequantization if and only if $K=\hbar\lb{N+\frac{n+1}{2}}$ for some $N\in\Z$.  Note that $\CP^n$ admits metaplectic-c prequantizations for all $n$.  This is an improvement over the Kostant-Souriau quantization scheme with half-form correction, because $\CP^n$ does not admit a metaplectic structure when $n$ is even.\footnote{Some additional detail:  if we consider a regular level set of the $n$-dimensional harmonic oscillator from Section \ref{subsec:SHO} corresponding to energy $K$, the symplectic reduction at this level is the symplectic manifold $(\CP^{n-1},Ki\varpi_{FS})$.  The fact that the reduced system admits a metaplectic-c prequantization when $K=\hbar\lb{N+\frac{n}{2}}$ can be checked directly using properties of $\CP^{n-1}$, as in \cite{rr1}, or it can be seen immediately by applying Theorem \ref{thm:MpcRed} to the $n$-dimensional harmonic oscillator.}

As a concrete example, we consider $\CP^2$, which does not admit a metaplectic structure.  For any $K$ of the form $K=\hbar\lb{N+\frac{3}{2}}$, $N\in\Z$, the symplectic manifold $(\CP^2,\omega_{FS})$ does not admit a prequantization line bundle either, since $\omega_{FS}$ is not integral.  However, it does admit a metaplectic-c prequantization for any such $K$.  We choose $K=\frac{3}{2}\hbar$.

Define an action of $T^2$ on $\C^3$ such that if $\tau=(\tau_1,\tau_2)\in T^2$, then $$\tau\cdot z=(z_0,\tau_1^{k_1^1}\tau_2^{k_2^1}z_1,\tau_1^{k_1^2}\tau_2^{k_2^2}z_2),$$ where $k^j=(k_1^j,k_2^j)\in\Z^2$ for $j=1,2$, and $\cb{k^1,k^2}$ is an integer basis for $\Z^2$.  This action descends to an effective Hamiltonian action of $T^2$ on $(\CP^2,\omega_{FS})$.  A calculation establishes that the momentum map takes the form $$\Phi([z])=- \frac{3h}{2}\lb{k_1^1\frac{|z_1|^2}{|z|^2}+k_1^2\frac{|z_2|^2}{|z|^2},k_2^1\frac{|z_1|^2}{|z|^2}+k_2^2\frac{|z_2|^2}{|z|^2}}+h(C_1,C_2),\ \ \forall [z]\in\CP^2,$$ where $(C_1,C_2)\in\R^{2*}$ is a constant.

The fixed points of the action are $Z_0=[1,0,0]$, $Z_1=[0,1,0]$ and $Z_2=[0,0,1]$.  It suffices to find a value of $(C_1,C_2)$ such that the value of $\Phi$ at one of these points satisfies the equivariance condition in Theorem \ref{thm:Equivar}.  For example, at $Z_0$, we have $$\frac{1}{h}\Phi(Z_0)=(C_1,C_2).$$  The equivariance condition at this point is $$\frac{1}{h}\Phi(Z_0)-\frac{1}{2}w_{Z_0*}=(C_1,C_2)-\frac{1}{2}w_{Z_0*}\in\Z^{2*},$$ where $w_{Z_0}:T^k\rightarrow U(1)$ is defined in terms of the action of $T^2$ on the tangent space $T_{Z_0}\CP^2$.  We can satisfy the equivariance condition by taking $(C_1,C_2)=\frac{1}{2}w_{Z_0*}$.  It remains to calculate $\frac{1}{2}w_{Z_0*}$.

On the open set $\cb{[z]\in\CP^2:z_0\neq 0}$, use local coordinates $(\zeta_1,\zeta_2)=\lb{\frac{z_1}{z_0},\frac{z_2}{z_0}}$ for $\CP^2$.  Let $\xi=\sum_{j=1}^2a_j\xi_j\in\mt$ be arbitrary, and let $\tau=\exp(\xi)=(e^{2\pi ia_1},e^{2\pi ia_2})$.  Then $\tau$ acts on the point $(\zeta_1,\zeta_2)$ by $$\tau\cdot(\zeta_1,\zeta_2)=(e^{2\pi i(a_1k^1_1+a_2k^1_2)}\zeta_1,e^{2\pi i(a_1k^2_1+a_2k^2_2)}\zeta_2).$$  
Identify $T_{Z_0}\CP^2$ with $\C^2$ in the natural way.  Then the complex matrix corresponding to the action of $\tau$ on $T_{Z_0}\CP^2$ is $$\mbox{diag}\lb{\exp\ls{2\pi i(a_1k_1^1+a_2k_2^1)},\exp\ls{2\pi i(a_1k_1^2+a_2k_2^2)}}.$$  Therefore $$w_{Z_0}(\tau)=\exp\ls{2\pi i(k_1^1+k_1^2)a_1+2\pi i(k_2^1+k_2^2)a_2},$$ which implies that $$\frac{1}{2}w_{Z_0*}=\frac{1}{2}\lb{k_1^1+k_1^2,k_2^1+k_2^2}\in\R^{2*}.$$  Hence we set $(C_1,C_2)=\frac{1}{2}\lb{k_1^1+k_1^2,k_2^1+k_2^2}$.

We can verify that this choice of $(C_1,C_2)$ also satisfies the equivariance condition over $Z_1$.  At $Z_1$, we have   $$\frac{1}{h}\Phi(Z_1)=-\frac{3}{2}\lb{k_1^1,k_2^1}+\frac{1}{2}\lb{k_1^1+k_1^2,k_2^1+k_2^2}.$$  We need to calculate $\frac{1}{2}w_{Z_1*}$, where $w_{Z_1}$ is defined in terms of the $T^2$ action on the tangent space $T_{Z_1}\CP^2$.  On the open set $\cb{[z]\in\CP^2:z_1\neq 0}$, use local coordinates $(\zeta_0,\zeta_2)=\lb{\frac{z_0}{z_1},\frac{z_2}{z_1}}$.  Then $\tau=(e^{2\pi ia_1},e^{2\pi ia_2})$ acts by $$\tau\cdot(\zeta_0,\zeta_2)=(e^{2\pi i(-a_1k_1^1-a_2k_2^1)}\zeta_0,e^{2\pi i(a_1(k_1^2-k_1^1)+a_2(k_2^2-k_2^1))}\zeta_2).$$  An identical calculation to the one performed at $Z_0$ yields $$\frac{1}{2}w_{Z_1*}=\frac{1}{2}(k_1^2-2k_1^1,k_2^2-2k_2^1).$$  Now, 
\begin{align*}
\frac{1}{h}\Phi(Z_1)-\frac{1}{2}w_{Z_1*}=&-\frac{3}{2}\lb{k_1^1,k_2^1}+\frac{1}{2}\lb{k_1^1+k_1^2,k_2^1+k_2^2}-\frac{1}{2}(k_1^2-2k_1^1,k_2^2-2k_2^1)\\
=&(0,0)\in\Z^{2*},
\end{align*}
as needed.  One can similarly check $Z_2$.

The image of the momentum map in $\R^{2*}$ (scaled by $\frac{1}{h}$, for simplicity) is the triangle with vertices 
\begin{align*}
\frac{1}{h}\Phi(Z_0)=&\frac{1}{2}\lb{k_1^1+k_1^2,k_2^1+k_2^2},&\frac{1}{h}\Phi(Z_1)=&\frac{1}{2}\lb{k_1^2-2k_1^1,k_2^2-2k_2^1},\\
\frac{1}{h}\Phi(Z_2)=&\frac{1}{2}\lb{k_1^1-2k_1^2,k_2^1-2k_2^2}.
\end{align*}
The quantized energy levels correspond to the integer lattice points lying strictly in the interior of the triangle.  

As a particularly simple example, consider $k^1=(1,0)$ and $k^2=(0,1)$.  The vertices of the image of the momentum map are $$\frac{1}{h}\Phi(Z_0)=\lb{\frac{1}{2},\frac{1}{2}},\ \ \frac{1}{h}\Phi(Z_1)=\lb{-1,\frac{1}{2}},\ \ \frac{1}{h}\Phi(Z_2)=\lb{\frac{1}{2},-1}.$$  There is exactly one integer lattice point in the interior of the triangle, namely $(0,0)$.  More generally, if we let $K=\hbar\lb{N+\frac{3}{2}}$ for an arbitrary $N\in\Z$, $N\geq 0$, the vertices are $$\frac{1}{h}\Phi(Z_0)=\lb{\frac{1}{2},\frac{1}{2}},\ \ \frac{1}{h}\Phi(Z_1)=\lb{-N-1,\frac{1}{2}},\ \ \frac{1}{h}\Phi(Z_2)=\lb{\frac{1}{2},-N-1},$$ and the number of quantized energy levels is $\frac{(N+2)(N+1)}{2}={{N+2}\choose{N}}$.  If $N=-1$, the system is metaplectic-c prequantizable and the vertices take the form given above, but there are no quantized energy levels.

The symplectic reduction of the three-dimensional harmonic oscillator at the quantized energy level $K=\hbar\lb{N+\frac{3}{2}}$ is exactly the symplectic manifold $(\CP^2,Ki\varpi_{FS})$.  We recognize the value ${{N+2}\choose{N}}$ from the quantum mechanical calculation as the multiplicity of the $N$th quantized energy level for $N\geq 0$.  This calculation does not yield a quantized energy level corresponding to $N=-1$, and indeed this regular value has multiplicity zero by the above interpretation.  

This last example and that in Section \ref{subsec:SHO} are different facets of the same system.  We will treat the relationships between them in greater generality in a subsequent paper concerning equivariant metaplectic-c prequantizations and quantized energy levels for toric manifolds (in preparation).

\end{document}